\documentclass{amsart}
\usepackage{amsmath,amsthm}
\usepackage{amsfonts,amssymb}
\usepackage{enumerate}
\usepackage{graphicx}

\newtheorem{lemma}{Lemma}

\newtheorem{corollary}[lemma]{Corollary}
\newtheorem{theorem}[lemma]{Theorem}
\newtheorem{remark}{Remark}

\begin{document}

\title{An inequality for Jacobi polynomials: a complement to Finite Increment Theorem}

\author{Geno Nikolov}
\address{Faculty of Mathematics and Informatics, Sofia University
"St. Kliment Ohridski", 5 James Bourchier Blvd., 1164 Sofia,
Bulgaria} \email{geno@fmi.uni-sofia.bg}

\begin{abstract}
Let $P=P_n^{(\alpha,\beta)}$ be the $n$-th degree Jacobi polynomial,
which is orthogonal in $[-1,1]$ with respect to the weight function
$(1-x)^{\alpha}(1+x)^{\beta}$, $\alpha,\beta>-1$. For parameters
$(\alpha,\beta)$ satisfying either $\alpha\geq\beta\geq 1/2$ or
$\alpha\geq 1/2$, $\beta=-1/2$, we prove the inequality
$$
P(1)-P(x)\geq P^{\prime}(x)\,(1-x),\quad x\in [0,1],
$$
which may be viewed as a complement to Finite Increment Theorem for
Jacobi polynomials.
  \vspace{5pt}
  \newline
  \textbf{Keywords: }{Finite Increment Theorem, Chebyshev polynomials,
  Jacobi polynomials, nonnegative expansion.}

  \vspace{5pt}
  \noindent
  \textbf{2020 Mathematics Subject Classification:} 41A17.
\end{abstract}

\maketitle \pagestyle{myheadings} \markboth{G. Nikolov}
         {Complement to Finite Increment Theorem}


\section{Introduction and statement of the result}
In \cite{GN2020} we investigated two trigonometric inequalities
proposed by M. S. Robertson \cite{MR1945} and R. Askey and G. Gasper
\cite{AG1976}, respectively (see (1.29) and (8.17) in
\cite{RA1975}). These inequalities were restated in \cite{GN2020} in
terms of the Chebyshev polynomial of the first kind
\[
T_n(x)=\cos (n\theta), \quad x=\cos\theta\,,
\]
and the following refinement was proved:
\medskip

\noindent {\bf Theorem A.} \textit{Let $n\in\mathbb{N}$, $n\geq 4$.
Then}
\begin{equation}\label{e1.1}
T_n(x)+2-\frac{x+2}{n^2}\,T_n^{\prime}(x)-
a(n)\frac{1-x}{n^2}\,\big(n^2-T_n^{\prime}(x)\big)\geq 0,\quad x\in
[-1,1],
\end{equation}
\textit{where}
\begin{equation}\label{e1.2}
 a(n)=\begin{cases}\frac{1}{1+\cos\frac{\pi}{n}}\,,&\textit{ if }\ n
 \ \textit{ is even},\\
\frac{1}{1+\cos\frac{2\pi}{n}}\,,&\textit{ if }\ n \
 \textit{ is odd.}
 \end{cases}
\end{equation}
\textit{The constant} $a(n)$ \textit{ in } \eqref{e1.2} \textit{ is
the best possible in the sense that  \eqref{e1.1} fails for any
larger constant. The equality in \eqref{e1.1} is attained only at
$x=1$ and $x=-\cos\frac{\pi}{n}$ if $n$ is even, and at $x=\pm 1$
and $x=-\cos\frac{2\pi}{n}$ if $n$ is odd.}
\medskip\noindent

The inequalities of Robertson and of Askey and Gasper are obtained
with choosing in \eqref{e1.1} $a(n)=0$ and $a(n)=1/2$, respectively.
We note that Theorem~A is true also for $n=1,\,2,\,3$, in which
cases the inequality in \eqref{e1.1} becomes an identity.

For $x\in [0,1]$ the inequality in \eqref{e1.1} holds true with a
larger constant $a(n)$, namely with $a=1$, yielding (cf.
\cite[Theorem~2]{GN2020})

\medskip\noindent
{\bf Theorem B.} \textit{For every $n\in\mathbb{N}$, $n\geq 3$, the
following inequality holds true:}
\begin{equation}\label{e1.3}
T_n^{\prime}(1)-T_n^{\prime}(x)\geq (1-x)\,T_n^{\prime\prime}(x),
\quad x\in [0,1].
\end{equation}
\textit{The equality in \eqref{e1.3} occurs only for $x=1$ and, if
$n\equiv 2\, (\!\!\!\mod 4)$, for $\,x=0$.}
\medskip

In the cases $n=0,\,1,\,2$, inequality \eqref{e1.3} becomes an
identity. The Chebyshev polynomials of the second kind $\{U_m\}$ are
related to the Chebyshev polynomial of the first kind by
\begin{equation}\label{e1.4}
U_{m}(x)=\frac{\sin
(m+1)\theta}{\sin\theta}=\frac{T_{m+1}^{\prime}(x)}{m+1}\,, \quad
x=\cos\theta.
\end{equation}
Therefore, inequality \eqref{e1.3} can be reformulated as
\begin{equation}\label{e1.5}
U_n(1)-U_n(x)\geq (1-x)\,U_n^{\prime}(x), \quad n\in\mathbb{N}_0,\
x\in [0,1].
\end{equation}

According to the Finite Increment Theorem, $\,U_n(1)-U_n(x)=
(1-x)\,U_n^{\prime}(\xi)$ for some $\xi\in (x,1)$, hence inequality
\eqref{e1.3} may be viewed as a complement to this theorem applied
to $U_n$.

The Chebyshev polynomials of the first and the second kind are
ultraspherical polynomials orthogonal in $(-1,1)$ with respect to
the weight functions $1/\sqrt{1-x^2}$ and $\sqrt{1-x^2}$,
respectively. On the other hand, the ultraspherical (called also
Gegenbauer) polynomials belong to the family of Jacobi polynomials
$\{P_n^{(\alpha,\beta)}(x)\}_{n\in\mathbb{N}_0}$, which are
orthogonal in $[-1,1]$ with respect to the weight function
$w_{\alpha,\beta}(x)=(1-x)^{\alpha}(1+x)^{\beta}$,
$\alpha,\,\beta>-1$. Inequality \eqref{e1.5} says that for
$(\alpha,\beta)=(1/2,1/2)$ we have
\begin{equation}\label{e1.6}
P_n^{(\alpha,\beta)}(1)-P_n^{(\alpha,\beta)}(x)\geq (1-x)\,
\frac{d}{dx}\big\{P_n^{(\alpha,\beta)}(x)\big\}\,, \quad x\in
[0,1]\,.
\end{equation}
A rather natural question is: for which other parameters
$(\alpha,\beta)$ the inequality \eqref{e1.6} holds true ?

The aim of this note is to provide a partial answer of the above
question. Our result is
\begin{theorem}\label{t1.1}
The inequality \eqref{e1.6} holds true in the following cases:
\begin{itemize}
\item[(i)~~]
$\alpha\geq\frac{1}{2},\ \beta=-\frac{1}{2}\,$;
\item[(ii)~~]
$\alpha\geq\beta\geq\frac{1}{2}\,$.
\end{itemize}
Unless $\alpha=\beta=\frac{1}{2}$, for $n\geq 3$ the equality in
\eqref{e1.6} is attained only at $x=1$.
\end{theorem}

The special case $\alpha=\frac{1}{2}$, $\beta=-\frac{1}{2}$ of
Theorem~\ref{t1.1}(i)  comes down to a trigonometric inequality,
which to the best of our knowledge is new:
\begin{corollary}\label{c1.1}
For every $n\in\mathbb{N},\;n\geq 2$, the following inequality holds
true:
$$
2n+1-(n+2)\,\frac{\sin
n\theta}{\sin\theta}+(n-1)\,\frac{\sin(n+1)\theta}{\sin\theta}\geq
0\,,\quad \theta \in \big[0,\frac{\pi}{2}\big]\,.
$$
The equality is attained only for $\theta=0$.
\end{corollary}

The proofs are given in the next section. We provide a short proof
of Theorem~B for two reasons: firstly, for the readers convenience,
and secondly, because Theorem~B constitutes the main ingredient of
our arguments.

\section{Proofs}
\subsection{Proof of Theorem~B}
Denote by $\tau$ the largest zero of $T_n^{\prime\prime}$, then
$$
\cos\frac{2\pi}{n}<\tau<\cos\frac{\pi}{n}\,.
$$
$T_n^{\prime\prime}(x)$ increases monotonically in $[\tau,1]$,
therefore
$$
T_n^{\prime}(1)-T_n^{\prime}(x)=\int\limits_{x}^{1}T_n^{\prime\prime}(u)\,du
\geq (1-x)T_n^{\prime\prime}(x),\quad x\in[\tau,1],
$$
with the equality only when $x=1$. Now the assertion of Theorem~B
obviously follows when $n=3$, and we assume in what follows $n\geq
4$. Since the left-hand side of \eqref{e1.3} is positive in $[0,1)$
while $\,T_n^{\prime\prime}(x)<0$ for $x\in
\big[\cos\frac{2\pi}{n},\tau\big)$, it follows that \eqref{e1.3}
holds true with $">"$ sign when
$x\in\big[\cos\frac{2\pi}{n},\tau\big)$. Thus, we need to prove the
inequality in \eqref{e1.3} for $x\in [0,\cos\frac{2\pi}{n})$ and
$n\geq 5$. Let us set
\begin{equation}\label{e2.1}
\mathcal{L}(x)=y^{\prime}(1)-y^{\prime}(x)-(1-x)y^{\prime\prime}(x)\,,
\qquad y=T_n(x)\,.
\end{equation}
By using the differential equation for $y$ (see, e.g., \cite[Eqn.
(1.92)]{TR2020})
\begin{equation}\label{e2.2}
(1-x^2)y^{\prime\prime}(x)-x\,y^{\prime}(x)+n^2y(x)=0
\end{equation}
and the representations
$$
T_n(x)=\cos n\theta, \ T_n^{\prime}(x)=n\,\frac{\sin
n\theta}{\sin\theta},\quad x=\cos\theta,
$$
we find
$$
\mathcal{L}(0)=\begin{cases}
2n^2,& \text{ if } n\equiv 0\,(\!\!\!\!\mod 4),\\
n^2-n,& \text{ if } n\equiv 1\,(\!\!\!\!\mod 4),\\
0,& \text{ if } n\equiv 2\,(\!\!\!\!\mod 4),\\
n^2+n,& \text{ if } n\equiv 3\,(\!\!\!\!\mod 4).
\end{cases}
$$
Thus, $\mathcal{L}(0)\geq 0$ and the equality is attained only when
$n=4k+2,\ k\in\mathbb{N}$. It remains to prove that
$\mathcal{L}(x)>0$ for $x\in (0,\cos\frac{2\pi}{n})$, and it
suffices to verify this inequality at the points of local extremum
of $\mathcal{L}(x)$ therein. Since
$\mathcal{L}^{\prime}(x)=(1-x)y^{\prime\prime\prime}(x)$, the latter
points are the zeros of $y^{\prime\prime\prime}(x)$ in
$\,(0,\cos\frac{2\pi}{n})$. Let $t\in (0,\cos\frac{2\pi}{n})\,$ be a
zero of $y^{\prime\prime\prime}(x)$. Note that
$$
\cos\frac{2\pi}{n}=1-2\sin^2\frac{\pi}{n}<
1-2\Big(\frac{2}{\pi}\frac{\pi}{n}\Big)^2=1-\frac{8}{n^2}\,,
$$
therefore
$$
0<t<1-\frac{8}{n^2}\,.
$$
Differentiating \eqref{e2.2} at $x=t$ and using
$y^{\prime\prime\prime}(t)=0$, we express $y^{\prime\prime}(t)$
through $y^{\prime}(t)$ as
$$
y^{\prime\prime}(t)=\frac{n^2-1}{3t}\,y^{\prime}(t).
$$
By plugging this expression in equation \eqref{e2.2} with $x=t$, we
obtain
\begin{equation}\label{e2.3}
y^{\prime}(t)=-\frac{3n^2 t}{n^2-1-(n^2+2)t^2}\,y(t)\,,
\end{equation}
and consequently
\begin{equation}\label{e2.4}
y^{\prime\prime}(t)=-\frac{n^2(n^2-1)}{n^2-1-(n^2+2)t^2}\,y(t)\,.
\end{equation}
Now we replace $y^{\prime}(1)=n^2$ and the expressions from
\eqref{e2.3} and \eqref{e2.4} in $\mathcal{L}(t)$ to obtain
\begin{equation}\label{e2.5}
\mathcal{L}(t)=n^2\Big[1+\frac{n^2-1-(n^2-4)t}{n^2-1-(n^2+2)t^2}\,y(t)\Big]\,.
\end{equation}
The quotient in front of $y(t)$ in the last expression is positive.
Indeed, this is obvious for the numerator, and for the denominator
follows from
$$
t^2<t<1-\frac{8}{n^2}<1-\frac{3}{n^2+2}=\frac{n^2-1}{n^2+2}\,.
$$
Since $y(t)\geq -1$, we conclude from \eqref{e2.5} that
$$
\mathcal{L}(t)\geq
n^2\Big[1-\frac{n^2-1-(n^2-4)t}{n^2-1-(n^2+2)t^2}\Big]=
\frac{t\big[n^2-4-(n^2+2)t)\big]}{n^2-1-(n^2+2)t^2}\,.
$$
The numerator of the last quotient is positive, which is seen from
$$
t<1-\frac{8}{n^2}<1-\frac{6}{n^2+2}=\frac{n^2-4}{n^2+2}\,.
$$
Hence, $\mathcal{L}(t)>0$ and the proof of Theorem~B is
accomplished.
\subsection{Proof of Theorem~\ref{t1.1}}
We start with the special case
$\,(\alpha,\beta)=\big(\frac{1}{2},-\frac{1}{2}\big)$ of
Theorem~\ref{t1.1}(i). Apart from a constant multiplier, the Jacobi
polynomial $P_n^{(1/2,-1/2)}(x)$ coincides with the Chebyshev
polynomial of the fourth kind
\begin{equation}\label{e2.6}
W_n(x)=\frac{\sin\big(n+\frac{1}{2}\big)\theta}{\sin\frac{\theta}{2}},\qquad
x=\cos\theta\,.
\end{equation}
This representation and \eqref{e1.4} imply
\begin{equation}\label{e2.7}
U_{n}(x)+U_{n-1}(x)=\frac{\sin(n+1)\theta}{\sin\theta} +\frac{\sin
n\theta}{\sin\theta}=
\frac{\sin\big(n+\frac{1}{2}\big)\theta}{\sin\frac{\theta}{2}}=W_n(x)\,.
\end{equation}
Now using \eqref{e1.5}, we obtain for $x\in [0,1]$
\begin{equation*}
\begin{split}
W_n(1)-W_n(x)&=U_n(1)-U_n(x)+U_{n-1}(1)-U_{n-1}(x)\\
&\geq (1-x)U_n^{\prime}(x)+(1-x)U_{n-1}^{\prime}(x)\\
&=(1-x)W_n^{\prime}(x)\geq 0\,.
\end{split}
\end{equation*}
For $n\geq 3$ the equality occurs only for $x=1$. Indeed, according
to Theorem~B, the equality can only occur at $x=0$ or $x=1$, but it
is not possible for both $U_n(1)-U_n(x)-(1-x)U_n^{\prime}(x)$ and
$U_{n-1}(1)-U_{n-1}(x)-(1-x)U_{n-1}^{\prime}(x)$ to vanish at $x=0$.

We shall need the next lemma, which provides a positive expansion
property of Jacobi polynomials.

\begin{lemma}\label{l2.1}
The Jacobi polynomials obey the following properties:
\begin{itemize}
\item[(i)]
If $\alpha>\gamma>-1$, then
$$
P_{n}^{(\alpha,\beta)}=\sum_{k=0}^{n}c_{k,n}\,P_{k}^{(\gamma,\beta)}
\quad \text{with all }\ c_{k,n}> 0\,.
$$
\item[(ii)]
If $\gamma>\alpha>-1$, then
$$
P_{n}^{(\gamma,\gamma)}=\sum_{k=0}^{n}c_{k,n}\,P_{k}^{(\alpha,\alpha)}
\quad \text{with all }\ c_{k,n}\geq 0\,,
$$
and $\,c_{k,n}>0\,$ if $\,k\equiv n\ (\!\!\!\mod 2)$.
\end{itemize}
\end{lemma}
The explicit formulae for the coefficients $\{c_{k,n}\}$ are given
in \cite[Eqns. (7.33) and (7.34)]{RA1975} (another reference for
part (ii), which concerns the ultraspherical polynomials, is
\cite[pp. 95--96]{GS1975}).

The general case of Theorem~\ref{t1.1}(i) is a consequence from the
case $(\alpha,\beta)=(1/2,-1/2)$ and Lemma~\ref{l2.1}(i) with
$\gamma=\frac{1}{2}$ and $\beta=-\frac{1}{2}$. If
$\alpha>\frac{1}{2}$ and $x\in [0,1]$, then
\begin{equation*}
\begin{split}
&P_n^{(\alpha,-1/2)}(1)-P_n^{(\alpha,-1/2)}(x)
-(1-x)\,\frac{d}{dx}\big\{P_n^{(\alpha,-1/2)}(x)\big\}\\
&=\sum_{m=0}^{n}c_{m,n}\,\Big[P_m^{(1/2,-1/2)}(1)-P_m^{(1/2,-1/2)}(x)-
(1-x)\,\frac{d}{dx}\big\{P_m^{(1/2,-1/2)}(x)\big\}\Big]\geq 0\,.
\end{split}
\end{equation*}
For $n\geq 3$ the equality is attained only at $x=1$ since all
summands with indices $m\geq 3$ in the above sum are positive at
$x=0$, and so is the sum.
\smallskip

Theorem~\ref{t1.1}(ii) follows from Theorem~B and Lemma~\ref{l2.1}
in a similar manner. For $\gamma\geq\frac{1}{2}$ and  $x\in [0,1]$
Lemma~\ref{l2.1}(ii) and Theorem~B imply
\begin{equation*}
\begin{split}
P_n^{(\gamma,\gamma)}(1)&-P_n^{(\gamma,\gamma)}(x)
-(1-x)\,\frac{d}{dx}\big\{P_n^{(\gamma,\gamma)}(x)\big\}\\
&=\sum_{m=0}^{n}c_{m,n}\,\Big[P_m^{(1/2,1/2)}(1)-P_m^{(1/2,1/2)}(x)-
(1-x)\,\frac{d}{dx}\big\{P_m^{(1/2,1/2)}(x)\big\}\Big]\\
&\geq 0\,.
\end{split}
\end{equation*}
If $n\geq 3$, then the equality is attained only at $x=1$. Indeed,
on account of Theorem~B and Lemma~\ref{l2.1}(ii), depending on
whether $n$ is odd or even, the summand corresponding to $m=3$ or
$m=4$ in the last sum is positive at $x=0$. Therefore, the sum is
positive at $x=0$, too.

Invoking again Lemma~\ref{l2.1}(i), we conclude from the case just
considered that if $\alpha>\beta\geq \frac{1}{2}$ and $x\in [0,1]$,
then
\begin{equation*}
\begin{split}
P_n^{(\alpha,\beta)}(1)-&P_n^{(\alpha,\beta)}(x)
-(1-x)\,\frac{d}{dx}\big\{P_n^{(\alpha,\beta)}(x)\big\}\\
&=\sum_{m=0}^{n}c_{m,n}\,\Big[P_m^{(\beta,\beta)}(1)-P_m^{(\beta,\beta)}(x)-
(1-x)\,\frac{d}{dx}\big\{P_m^{(\beta,\beta)}(x)\big\}\Big]\geq 0
\end{split}
\end{equation*}
and the equality is attained only for $x=1$. Theorem~\ref{t1.1} is
proved.
\subsection{Proof of Corollary~\ref{c1.1}}
We show that Corollary~\ref{c1.1} is an alternative representation
of the inequality
\begin{equation}\label{e2.9}
W_n(1)-W_n(x)-(1-x)W_n^{\prime}(x)\geq 0,\quad x\in [0,1]\,.
\end{equation}
From \eqref{e2.6} we obtain
\begin{equation*}
\begin{split}
W_n^{\prime}(x)&=\frac{-1}{\sin\theta}\,\frac{d}{d\theta}\Big\{\frac{\sin
\big(n+\frac{1}{2}\big)\theta}{\sin\frac{\theta}{2}}\Big\}\\
&=\frac{1}{\sin\theta\sin^2\frac{\theta}{2}}
\Big[\frac{1}{2}\sin\Big(n+\frac{1}{2}\Big)\theta\cos\frac{\theta}{2}-
\Big(n+\frac{1}{2}\Big)\sin\frac{\theta}{2}\cos\Big(n+\frac{1}{2}\Big)\theta
\Big]\\
&=\frac{1}{2\sin\theta\sin^2\frac{\theta}{2}}\Big[
\frac{1}{2}\big[\sin(n+1)\theta+\sin(n\theta)\big]-
\Big(n+\frac{1}{2}\Big)\big[\sin(n+1)\theta-\sin(n\theta)\big]
\Big]\\
&=\frac{1}{2\sin^2\frac{\theta}{2}}\, \Big[
(n+1)\,\frac{\sin(n\theta)}{\sin\theta}
-n\,\frac{\sin(n+1)\theta}{\sin\theta}\Big]\,.
\end{split}
\end{equation*}
Consequently,
$$
(1-x)\,W_n^{\prime}(x)=(n+1)\,\frac{\sin(n\theta)}{\sin\theta}
-n\,\frac{\sin(n+1)\theta}{\sin\theta} \,.
$$
By replacing $W(x)$ from \eqref{e2.6} and plugging the above
expression in \eqref{e2.9}, we arrive at the equivalent
trigonometric inequality
\begin{equation}\label{e2.10}
2n+1-(n+2)\,\frac{\sin
n\theta}{\sin\theta}+(n-1)\,\frac{\sin(n+1)\theta}{\sin\theta}\geq
0\,,\quad \theta \in \big[0,\frac{\pi}{2}\big]\,,
\end{equation}
with the equality sign occurring only for $\theta=\pi/2$. The proof
of Corollary~\ref{c1.1} is complete. Let us point out to the
equivalent to Corollary~\ref{c1.1} statement given in terms of the
Chebyshev polynomials of the second kind:
\begin{corollary}\label{c2.1}
For every $n\in\mathbb{N},\;n\geq 2$, the following inequality holds
true:
$$
(n+2)\big[U_{n-1}(1)-U_{n-1}(x)\big]-(n-1)\big[U_{n}(1)-U_{n}(x)\big]\geq
0,\quad x\in [0,1]\,.
$$
The equality is attained only for $x=1$.
\end{corollary}
We conclude this note with the following comment.
\begin{remark}
A challenging task is to find domain $\mathcal{D}\subset
\mathbb{R}^2$ formed by the pairs $(\alpha,\beta)$ such that for all
$\,n\in\mathbb{N},\;n\geq 3$ the inequality
$$
P_n^{(\alpha,\beta)}(1)-P_n^{(\alpha,\beta)}(x)\geq (1-x)\,
\frac{d}{dx}\big\{P_n^{(\alpha,\beta)}(x)\big\}\,, \quad x\in
[0,1]\,,
$$
holds true. Some experiments carried out with \emph{Wolfram
Mathematica} indicate that the region provided by Theorem~\ref{t1.1}
is a proper subset of $\mathcal{D}$, e.g., if, say, $\beta>2$, then
there are pairs $(\alpha,\beta)\in \mathcal{D}$ with $\alpha<\beta$.
\end{remark}


\section{Acknowledgment}
This research is partially supported by the Bulgarian National
Research Fund under Contract~KP-06-N62/4


\end{document}